
\tolerance=10000
\raggedbottom

\baselineskip=15pt
\parskip=1\jot

\def\sk{\vskip 3\jot}

\def\heading#1{\vskip3\jot{\noindent\bf #1}}
\def\label#1{{\noindent\it #1}}
\def\QED{\hbox{\rlap{$\sqcap$}$\sqcup$}}


\def\ref#1;#2;#3;#4;#5.{\item{[#1]} #2,#3,{\it #4},#5.}
\def\refinbook#1;#2;#3;#4;#5;#6.{\item{[#1]} #2, #3, #4, {\it #5},#6.} 
\def\refbook#1;#2;#3;#4.{\item{[#1]} #2,{\it #3},#4.}


\def\({\bigl(}
\def\){\bigr)}


\def\al{\alpha}
\def\be{\beta}

\def\de{\delta}
\def\ep{\varepsilon}

\def\ph{\phi}

\def\Me{\Omega}


\def\calX{{\cal X}}


\def\bfD{{\bf D}}

\def\bfZ{{\bf Z}}

\def\abs#1{\vert#1\vert}

{
\pageno=0
\nopagenumbers
\rightline{\tt high.fault.rate.tex}
\vskip1in

\centerline{\bf Fault Tolerance in Cellular Automata at High Fault Rates}
\vskip0.5in

\centerline{Nicholas Pippenger}
\centerline{\tt njp@math.hmc.edu}
\centerline{Department of Mathematics}
\centerline{Harvey Mudd College}
\centerline{1250 Dartmouth Avenue}
\centerline{Claremont, CA 91711}
\vskip0.5in

\centerline{Mark McCann}
\centerline{\tt mmccann@cs.princeton.edu}
\centerline{Department of Computer Science}
\centerline{Princeton University}
\centerline{35 Olden Street}
\centerline{Princeton, NJ 08540}
\vskip0.5in

\noindent{\bf Abstract:}
A commonly used model for fault-tolerant computation is that of cellular automata.
The essential difficulty of fault-tolerant computation is present in the special case of
simply remembering a bit in the presence of faults, and that is the case we treat in this paper.
We are concerned with the degree (the number of neighboring cells on which the state
transition function depends) needed to achieve fault tolerance when the fault rate is high (nearly $1/2$).
We consider both the traditional transient fault model (where faults occur independently in time and space) and a recently introduced combined fault model which
also includes manufacturing faults (which occur independently in space, but which affect cells for all
time).
We also consider both a purely probabilistic fault model (in which the states of cells are perturbed
at exactly the fault rate) and an adversarial model (in which the occurrence of a fault gives control
of the state to an omniscient adversary).
We show that there are cellular automata that can tolerate 
a fault rate $1/2 - \xi$ (with $\xi>0$)  with degree $O\((1/\xi^2)\log(1/\xi)\)$,
even with adversarial combined faults.
The simplest such automata are
based on infinite regular trees, but our results also apply to other
structures (such as hyperbolic tessellations) that contain infinite regular trees.
We also obtain a lower bound of $\Me(1/\xi^2)$, 
even with purely probabilistic transient faults only.
\vfill\eject
}

\heading{1.  Introduction}

The theory of fault-tolerant computation has a history almost as old as that of fault-tolerant
communication.
The most widely used theoretical model for computation, the Turing machine, is unsuitable
for the study of fault-tolerant computation: it calls for leaving large amounts data unattended
on tapes for long periods of time; it seems unrealistic to assume that this data will not be corrupted
by failures (at a positive constant rate per tape cell and per time step), but for a Turing machine
(which can perform just one basic action per time step) there is no hope of keeping up with the
failures that occur in the absence of such an assumption.
The first study of fault-tolerance in a suitable computational model was undertaken by
von Neumann [N1], who used the model of combinational circuits.
These circuits are built from gates interconnected by wires in an acyclic fashion, so that
information flows unidirectionally from input terminals to output terminals, and each gate acts just once in any given computation by a circuit.

Von Neumann's most fundamental result is this:
for every error probability $\de>0$, there exists a failure probability $\ep>0$ such that
for every circuit  that performs some computation in the absence of failures, there exists
another circuit (in general,  deeper by a constant factor)
that performs the same computation, with error probability at most $\de$,
even if each gate in the new circuit suffers a fault independently with probability $\ep$.
A key feature of this result is that although $\ep$ depends on $\de$, it does not depend
on the size of the original circuit or on the complexity of the computation it performs
(though it does depend on the choice of the set of types of gates that are used in constructing the circuits).
In stating this result, we have used a convention that will be employed throughout this paper:
the term ``failure'' refers to a situation in which a component (such as a gate) does not perform its proper function; the term ``error'' refers to a situation in which some value of some variable (such as the signal carried on a wire) differs from the value it would have in the absence of any failures.

In this paper we shall deal exclusively with Boolean or binary information, for which signals or states
can assume only two possible values.
Thus a failure can occur in only one way: the value of a function is replaced by its complementary value.

As described above, von Neumann considered the case in which each gate failed independently with
some fixed probability $\ep$.
He also mentioned, however, the desirability of considering another failure model,
one in which the failures at the various gates are arbitrary, but subject to the constraint that they are stochastically dominated by independent random events that each occur with probability $\ep$.
We shall refer to these events as ``faults''.
There are two ways of looking at this new failure model.
One is that an adversary, who knows the inputs to the circuit, chooses the joint probability distribution for all the failures,  subject to the stochastic-domination constraint described above.
An alternative, which will be employed in this paper, is that the faults occur independently, and then
an adversary, who knows both the inputs to the circuit and the locations of the faults, decides
which failures will occur, subject to the constraint that a failure can occur at a given gate only
if a fault occurs at that gate.

We thus will deal with two failure models: the {\it purely probabilistic\/} failure model, in which
a failure occurs if and only if a fault occurs, and the {\it adversarial\/} fault model described
in the preceding paragraph.
There are several reasons for considering the adversarial model.
One is that it prevents faults from providing a benefit to a circuit.
(In the purely probabilistic model, the faults provide a source of random events to the circuit.
Since there are many known examples of randomized algorithms that outperform their 
best known deterministic
counterparts, the possibility exists that, with the purely probabilistic model, fault-tolerant circuits
might be smaller than any non-fault-tolerant circuits performing the same computations.
This would be an interesting phenomenon, but it is not the one we want to study under the name
``fault-tolerance''.)
A second reason for using the adversarial model is that it may be technically convenient.
(Negative results that hold for the purely probabilistic model also hold for the adversarial model,
since the adversary can always cause a failure whenever there is a fault.
But some negative results (though not the one in this this paper) appear to require an adversarial
model.
Surprisingly, the adversarial model may also be more convenient for proving positive results.
In the purely probabilistic model, for example, it may not be possible to construct a circuit
that ``simulates'' a gate (because the error probability of a circuit may depend on the 
values of the inputs to the circuit,
whereas the failure probability of a gate should not); this makes it difficult to prove ``change of basis''
results (see Pippenger [P2]) that are easily proved for the adversarial model.)
Finally, the adversarial model may be preferable simply because it is more realistic (or at least less unrealistic) in a given situation.
(Failures will not in practice occur with exactly equal probabilities and complete independence.
The adversarial model provides a measure of insurance against departures from these assumptions.)

As mentioned above, each gate in a circuit acts just once during any given computation.
Thus the combinational circuit model is unsuitable for the study of temporal effects, stemming from the
independence (or lack thereof) of faults in a given component at different times.
One possible model for the study of these effects is that of sequential circuits,
which may contain flip-flops as memory elements, and in which the assumption that there are no 
cycles is weakened to the assumption that there are no cycles that do not pass through flip-flops.
The input-output conventions used for sequential circuits are sometimes different from those used for
combinational circuits: the circuit may have no input or output terminals; rather the initial states of all
the flip-flops may be regarded as the ``input'', and their states at some later time as the ``output''
(see Kuznetsov [K], for example).

More commonly, however, temporal effects are studied through the model of cellular automata,
which were introduced by Ulam [U] and von Neumann [N2].
A {\it cellular automaton\/} is based on a directed graph (called the {\it lattice\/}).
Associated with each vertex $v$ in the graph is a {\it cell}, which is characterized by (1) a {\it state set\/}
$\calX_v$, (2) a {\it transition function\/} $\ph_v$, and (3) a one-to-one correspondence between the
argument positions of the transition function and the edges directed out of $v$: if edges are directed
from $v$ to $w_1, w_2, \ldots, w_k$, then the transition function is a map 
$\ph_v : \calX_{w_1}\times \calX_{w_2}\times\cdots\times \calX_{w_k}\to \calX_v$.

A {\it configuration\/} $x$ of a cellular automaton is an assignment of a state $x_v \in \calX_v$ 
to each cell $v$.
The configuration of a cellular automaton evolves in time (assumed to take non-negative integer values) in the following way.
The initial configuration $x(0)$ is assumed to be given.
Given the states $x_v(t)$ of the cells at time $t\ge 0$, their states at time $t+1$ are determined
by applying their transition functions to the states of their neighbors:
$x_v(t+1) = \ph_v\( x_{w_1}(t), x_{w_2}(t), \ldots, x_{w_k}(t) \)$.
We shall adopt the convention that the initial configuration is the ``input'' to the automaton, and that 
its configuration a some later time is its ``output''.

In this paper, we deal exclusively with {\it binary\/} automata, for which each cell has just two states:
$\calX_v = \{0,1\}$ for all $v$.
Thus each transition function is a Boolean function of the appropriate number of arguments.

Let us consider some examples at this point.
Our first example is Conway's ``Game of Life'' (see Berlekamp, Conway and Guy [B]).
Take as the lattice the graph having as vertices the points of the plane with integer coordinates
(that is, the elements of $\bfZ\times\bfZ$), and edges directed from each vertex to itself and to each of its eight nearest neighbors in the plane.
The transition function for each cell is the following: the next state of a cell is $1$ if its current state is $0$
and exactly three of its eight neighbors are in state $1$, or if its current state is $1$ and
either two or three of its neighbors are in state $1$.
The automorphism group of this automaton is generated by the translations $\bfZ\times\bfZ$ together 
with the dihedral group $\bfD_4$ of symmetries of the square.
A feature of this automaton is that it is computationally universal: started in an appropriate initial
configuration, it will simulate an arbitrary Turing machine.

Thus far we have dealt with deterministic cellular automata.
To discuss fault tolerance, we must consider automata with probabilistically occurring faults.
Probabilistic cellular automata were first considered by Stavskaya and Pyatetski\u{\i}-Shapiro [S],
and with an adversarial fault model by Toom [T2, T3].

As a second example, we consider ``Toom's Rule'' (see Toom [T1, T2]).
Take the lattice to be the graph with the same vertices as in Conway's Game, but with edges from a 
vertex to itself and its ``northern'' and ``eastern''  neighbors, and take the transition function to be majority voting:  the next state of a cell is $1$ if and only if at least two of its three neighbors are in state $1$.
The automorphism group of this automaton is generated by the translations and the 
reflection about the main diagonal that exchanges the two coordinates;
it is {\it not\/} invariant under any rotations.
It is not hard to see that any initial configuration in which only finitely many cells are in state $1$
will eventually be driven to the all-$0$s configuration by iteration of the transition function,
and any configuration with only finitely many $0$s will be driven to all-$1$s.
Toom showed that it also has a more subtle property:   it can ``remember a bit'' forever, even in the presence of adversarial faults occurring at a sufficiently small rate.
That is, for every $\de>0$, there exists an $\ep>0$ such that, if all states are initially $a\in\{0,1\}$,
and if adversarial faults occur at rate $\ep$ (that is, the adversary is given control of the value of the transition function at each cell and each time independently with probability $\ep$), then
the probability that any given cell is in error (that is, is in state $1-a$) at any given time is at most $\de$.
This property of remembering a bit is all that is needed to achieve fault-tolerant computation:
by considering a cellular automaton based on a four-dimensional lattice, applying Coway's Game
in two of the dimensions and Toom's Rule in the other two, we obtain an automaton that simulates
an arbitrary Turing machine, with the state of each cell having arbitrarily small error probability when the fault rate is sufficiently small.

In the case of purely probabilistic faults, Toom's result amounts to showing that the stochastic process
associated with the cellular automaton and its probabilistic failures is non-ergodic, and that the all-$0$s and all-$1$s configurations lie in the basins of attraction of distinct invariant distributions on the configurations.
In the adversarial case, the presence of an adversary that can see into the future prevents the faulty
automaton from being considered as an autonomous stochastic process, but a special property of the 
transition function allows a reduction to an autonomous situation.

A Boolean function $\ph : \{0,1\}^k \to \{0,1\}$ is said to be {\it monotone\/} if increasing the value of an
argument from $0$ to $1$ cannot decrease the value of the function from $1$ to $0$: if
$x_1\le y_1, x_2\le y_2, \ldots, x_k\le y_k$, then $\ph(x_1, x_2, \ldots, x_k)\le \ph(y_1, y_2, \ldots, y_k)$.
Suppose a cellular automaton is started in the all-$0$s configuration, and that its transition function is monotone.
Then an adversary who is trying to maximize the probability that a particular cell is in state $1$ at a particular time has a clear optimal strategy: seize any opportunity to make the state of a cell $1$
(but decline any opportunity to make the state of a cell $0$),
for by monotonicity doing so cannot foreclose any future opportunities.
Similarly, if the automaton is started in the all-$1$s configuration, the adversary should seize any opportunity to make the state of a cell $0$.
The existence of these optimal ``greedy'' strategies means that for cellular automata with monotone transition functions, the analysis of adversarial faults can be reduced to the analysis of two stochastic processes: one in which all-$0$s is the initial configuration and a fault forces a state to $1$, and the
other in which all-$1$s is the initial configuration and a fault forces a state to $0$.

The majority voting function that is used in Toom's Rule has a property that further simplifies analysis:
it is self-dual.
A Boolean function $\ph : \{0,1\}^k \to \{0,1\}$ is said to be {\it self-dual\/} if it is invariant under exchanging
the roles of $0$ and $1$: $\ph(1-x_1, 1-x_2, \ldots, 1-x_k) = 1 - \ph(x_1, x_2, \ldots, x_k)$.
For a transition function that is self-dual as well as monotone, only one of the two stochastic processes described above needs to be considered.
Majority voting with any odd number of votes is both monotone and self-dual, and thus it plays
an important role in the construction of fault tolerant systems.

As described above, faults (either purely probabilistic or adversarial) in cellular automata are assumed
to occur independently both from time to time and from cell to cell.
This assumption is appropriate for studying {\it transient\/} faults, which affect the state of a cell
but do not impair its ability to function correctly in the future.
In practice, however, some types of faults {\it do\/} affect the functioning of cells.
To deal with these faults, McCann [M] has introduced a fault model that incorporates both
transient faults (as described above) and {\it manufacturing\/} faults, which are assumed to 
occur independently from cell to cell, but which when they occur at a cell give control of that
cell's state to the adversary for all time.
In the {\it combined\/} fault model, transient faults are assumed to occur (independently in time and space) at rate $\al>0$, and manufacturing faults are assumed to occur (independently in space)
at rate $\be>0$.
In the analysis, usually only the combined fault rate $\ep = 1 - (1-\al)(1-\be)$ (the probability that
a particular cell is subject to either a transient or a manufacturing fault at a particular time) is 
important.

McCann [M] has shown that Toom's Rule is not tolerant of combined faults (no matter how small
the fault rate) and indeed that no monotone binary cellular automaton based on the two-dimensional lattice $\bfZ\times\bfZ$ can tolerate combined faults.
He has also shown that a simple three-dimensional analog of Toom's Rule {\it  is\/} tolerant of combined faults.
This difference between two and three dimensions is significant because G\'{a}cs [G] has argued
that while two dimensional arrays of components are physically realistic, three dimensional ones are not, since they would require cubic amounts of power and heat to be transported through a 
boundary of quadratic area.

In this paper, we shall address the question of how large the degree (the number of neighbors
on which the transition function of a cell depends) must be for the automaton to tolerate faults at a fault rate
very close to $1/2$, that is for $\ep = 1/2 - \xi$ for some small $\xi>0$.
We shall obtain both upper and lower bounds to the degree.
The upper bounds will be obtained for highly structured automata, and under the hypotheses least favorable to fault tolerance: adversarial combined faults.
The graphs on which the automata are based will be undirected (an undirected edge comprises
two oppositely directed edges), regular and planar: they will be tessellations of the hyperbolic plane,
and thus they will have very large automophism groups.
Since these graphs are planar, these results contrast with McCann's negative result for the Euclidean
plane mentioned in the preceding paragraph.
The transition functions will be full majority voting: the new state of a cell is given by a majority
vote among the states of all its neighbors in the graph, including its own state if the number of neighbors is even.
Thus the transition functions will be both monotone and self-dual.
In Section 2, we shall describe automata meeting these criteria and having degree 
$O\((1/\xi^2)\log(1/\xi)\)$.

The lower bounds will be obtained under hypotheses most favorable to fault tolerance:
the graphs underlying the automata need not be planar or regular, and need not have any non-trivial
automorphisms, the transition functions need not be monotone or self-dual, and the automata need tolerate only purely probabilistic transient faults.
In Section 3, we shall show that even under these weak assumptions, the degee must be
$\Me(1/\xi^2)$.
\vfill\eject

\heading{2. Positive Results}

In this section, we shall construct cellular automata using full majority voting that tolerate
adversarial combined faults with fault rate $\ep = 1/2 - \xi$, error probability at most $\de = 1/2 - \xi/2$ and degree $O\((1/\xi^2)\log(1/\xi)\)$.
These automata will be based on highly symmetric undirected graphs (though the results will
also apply to unsymmetrical graphs), but the key to their fault tolerance will be a proposition
concerning automata based on directed trees (in which all edges are directed away from a root).

For $a\in\{0,1\}$, we define the {\it $a$-threshold\/} of a monotone Boolean function $\ph$ to be
the minimum number number of arguments of $\ph$  that, when set to $a$, force the value of $\ph$
to be $a$.
We define the {\it threshold\/} of $\ph$ to be the minimum of its $0$-threshold and its $1$-threshold.

\label{Proposition 2.1:}
Consider a monotone cellular automaton based on a directed tree $T$.
For every cell $v$, let $d(v)$ denote the out-degree of $v$,
and let $h(v)$ denote the threshold of the transition function $\ph_v$.
Suppose that for some $0<\xi<1/2$ and integer $m\ge 0$ we have an integer $d$ satisfying
$$d\ge m + {2\over \xi^2}\log{2^{m+1}\over \xi}.$$
Then if $d(v)\ge d$ and $h(v)\ge \(d(v)-m\)\big/2$ for all cells $v$, the automaton 
will tolerate adversarial combined faults with fault rate $\ep = 1/2 - \xi$ and error probability at most 
$1/2 - \xi/2$.

\label{Proof:}
Suppose, without loss of generality, that the value $0$ is to be remembered, so that all cells are initially in state $0$, and cells are in error when and only when they are in state $1$.
For $t\ge 0$, define $P_t$ to be the supremum over all cells $v$ of the probability that
$v$ is in error at time $t$.
We shall prove by induction on $t$ that if the fault rate is at most $\ep = 1/2 - \xi$, then
$P_t \le 1/2 - \xi/2$.
The base case is $P_0 = 0\le 1/2 - \xi/2$.

We now assume the bound $P_t \le 1/2 - \xi/2$, and prove $P_{t+1} \le 1/2 - \xi/2$.
If cell $v$ is in error at time $t+1$, then either (1) a fault occurs at $v$ at time $t+1$, or (2)
at least $h(v)$ of $v$'s children must have been in error a time $t$.
For $w$ a child of $v$,
let $E_w$ denote the event that cell $w$ is in error at time $t$.
Since there are no directed paths between distinct children of $v$,
the $d(v)$ events $E_w$ are independent.
Furthermore, since $\Pr[E_w]\le 1/2 - \xi/2$ by the inductive hypothesis, the $d(v)$ events
are stochastically dominated by $d(v)$ events that occur independently with probability
{\it exactly\/} $P_t$.
Thus we have
$$P_{t+1} \le \ep + \sum_{h(v)\le k\le d(v)} {d(v)\choose k} P_t^k (1-P_t)^{d(v)-k}.$$
Since $P_t < 1 - P_t$ and $\sum_k {d(v)\choose k} = 2^{d(v)}$, we have
$$\eqalign{
P_{t+1} &\le \ep + P_t^{h(v)} (1-P_t)^{d(v)-h(v)}\,\sum_{h(v)\le k\le d(v)} {d(v)\choose k} \cr
&\le \ep + P_t^{h(v)} (1-P_t)^{d(v)-h(v)}\,2^{d(v)} \cr
&= \ep + \(P_t /  (1-P_t)\)^{h(v)}\,\(2(1-P_t)\)^{d(v)}. \cr
}$$
Since $2(1-P_t) > 1$ and $d(v)\le 2h(v) + m$, we have
$$\eqalign{
P_{t+1} &\le \ep + \(P_t /  (1-P_t)\)^{h(v)}\,\(2(1-P_t)\)^{2h(v) + m} \cr
&= \(2(1-P_t)\)^m \, \(4P_t(1-P_t)\)^{h(v)}. \cr
}$$
Since $1-P_t < 1$, $4P_t(1-P_t) < 1$ and $h(v)\ge (d(v)-m)/2\ge (d-m)/2$, we have
$$P_{t+1} \le \ep + 2^m \, \(4P_t(1-P_t)\)^{(d-m)/2}.$$
We have $\ep = 1/2 - \xi$ and, by the inductive hypothesis, $P_t\le 1/2 - \xi/2$, so we obtain
$$P_{t+1} \le 1/2 - \xi + 2^m \, (1 - \xi^2)^{(d-m)/2}.$$
Thus to prove $P_{t+1} \le 1/2 - \xi/2$, it will suffice to show that
$$2^m \, (1 - \xi^2)^{(d-m)/2} \le \xi/2.$$
This inequality follows from the hypothesis of the proposition and the inequality $1-\xi^2 < \exp(-\xi^2)$:
$$\eqalign{
2^m \, (1 - \xi^2)^{(d-m)/2}
&\le 2^m \,\exp\left(-{\xi^2(d-m)\over 2}\right) \cr
&\le 2^m \,\exp\left(-\log{2^{m+1}\over \xi}\right) \cr
&= \xi/2. \cr
}$$
\QED

The following theorem extends the result of Proposition 2.1 to graphs that merely {\it contain\/}
a directed tree.

\label{Theorem 2.2:}
Consider a cellular automaton based on a graph $G$, with each vertex
having odd out-degree at least $s$, and the transition function for each cell being the majority function.
Suppose that it is possible to convert $G$ into a directed tree $T$ by deleting edges,
with at most $r$ of the edges directed out of any vertex being deleted.
Then if
$$s\ge 3r - 1 + {2\over \xi^2}\log{2^{2r}\over \xi},$$
the automaton will tolerate adversarial combined faults with fault rate $\ep = 1/2 - \xi$ 
and error probability $1/2 - \xi/2$.

\label{Proof:}
Suppose, without loss of generality, that the value $0$ is to be remembered, so that all cells are initially in state $0$, and cells are in error when and only when they are in state $1$.
Our strategy will be to delete edges from $G$ to convert it to $T$,
Whenever we delete an edge directed from a cell $v$ to a cell $w$, we will substitute the constant $1$
for the corresponding argument of $\ph_v$.
Since the constant $1$ stochastically dominates the actual state of $w$, an upper bound for the error 
probability in the tree automaton will also be an upper bound for the error probability in the original graph automaton.
To bound the error probability in the tree automaton, we estimate the out degrees of its vertices and
the thresholds of its transition functions.
These transition functions, being obtained from monotone functions by substitution
of constants for arguments, are themselves monotone, so
we may then apply Proposition 2.1.

Obviously each vertex $v$ of $T$ has out-degree  $d(v) \ge s - r$,
so the condition $d(v)\ge d$ of Proposition 2.1 will be fulfilled if we take $d = s - r$.
The transition function of the cell at $v$ in $G$ has threshold $\(d(v)+1\)\big/2$,
since it takes a majority of $d(v)$ votes.
The transition function of the cell at $v$ in $T$ therefore has threshold at least $\(d(v)+1\)\big/2 - r$.
Thus if we take $m = 2r - 1$, the condition $h(v)\ge\(d(v)-m\)\big/2$ of Proposition 2.1 will be fulfilled.
Finally, the condition
$$s\ge m + {2\over \xi^2}\log{2^{m+1}\over \xi}$$
of Proposition 2.1 will then be fulfilled by the hypothesis of the theorem.
\QED

In the following corollaries, we consider undirected graphs.
Each undirected edge will be regarded as two oppositely directed edges, and vertices
with even degree will be regarded as having a directed self-loop that represents
their inclusion in their own majority vote.

\label{Corollary 2.3:}
The regular $q$-ary tree with full majority voting tolerates adversarial combined faults with fault rate 
$\ep = 1/2 - \xi$ if $q$ is odd and
$$q \ge 2 + {2\over \xi^2}\log{4\over \xi},$$ 
or if $q$ is even and 
$$q \ge 4 + {2\over \xi^2}\log{16\over \xi}.$$ 

\label{Proof:}
Suppose first that $q$ is odd.
To convert the $q$-ary tree to a directed tree, is suffices to classify vertices into ``shells''
according to their distance (as measured by
the number of edges on a shortest path)
from an arbitrarily chosen  root, and to delete all edges that are directed from a farther vertex to a nearer one.
This amounts to deleting the edge from each child to its parent, and we may then apply Theorem 2.2 with  $s = q$ and $r = 1$.
If $q$ is even, we must include the self-loops to obtain a directed graph with out-degree $s = q + 1$.
To obtain a directed tree, we must delete the self-loop as well as the edge directed to the parent.
We then apply Theorem 2.2 with $r = 2$.
\QED

Regular trees of high degree have no cycles, but have ``expansion'', which manifests itself as a 
large ``isoperimetric constant'' (any finite set of vertices is adjacent to a proportional number of edges
that leave the set).
These trees are thus naturally imbedded in the hyperbolic plane.
That it is the expansion, and not the absence of cycles, that is the essential requirement for fault tolerance with majority voting is illustrated by examples based on regular hyperbolic tessellations (also known as ``honeycombs''), as described by Coxeter [C1].
In Coxeter's notation, $\{p,q\}$ (for $p\ge 3$ and $q\ge 3$ with $(p-2)(q-2) > 4$) denotes a tessellation of the hyperbolic plane in which $q$ $p$-gons meet at each vertex.
The automorphism groups of these tessellations are discussed by Coxeter and Moser [C2].
(For $(p-2)(q-2) = 4$, the notation $\{p,q\}$ denotes a regular tessellation of the Euclidean plane, and it is easy to see that
the corresponding cellular automata with majority voting are not tolerant of even purely probabilistic
transient faults.
For $(p-2)(q-2) < 4$,
$\{p,q\}$ denotes a regular tessellation of the sphere (that is, a Platonic solid), and of course the corresponding cellular automata, being finite, cannot remember a bit:
even with purely probabilistic transient failures, the stochastic process is ergodic, with each state
of each uniformly (but not independently) distributed in the invariant distribution on configurations.
The case $p=\infty$ corresponds to the $q$-ary tree.)

\label{Corollary 2.4:}
The cellular automaton using majority voting and based on
the tessellation $\{p,q\}$ with  $p\ge 3$,  tolerates combined adversarial faults with fault rate  
$\ep = 1/2 - \xi$ if $q$ if is odd and
$$q \ge 14 + {2\over \xi^2}\log{1024\over \xi},$$ 
or if $q$ is even and 
$$q \ge 16 + {2\over \xi^2}\log{4096\over \xi}.$$

\label{Proof:}
Suppose first that $q$ is odd.
Choose an arbitrary root vertex $v$, and classify vertices into shells
according to their distance  from $v$.
We count the number of directed edges that might have to be deleted to obtain a directed tree.
Consider a vertex $w$ in shell $n\ge 1$.
There can be at most two edges directed from $w$ to vertices in shell $n-1$ (the ``parents'' of $w$)
and at most two edges directed from $w$ to other vertices in shell $n$ (the ``siblings'' or ``cousins''
of $w$).
Finally, of the ``children'' of $w$ (the vertices in shell $n+1$ to which edges from $w$ are directed), we might have to exclude one, to ensure  that the remaining children of $w$ are disjoint from those of other
vertices in shell $w$.
In this way we delete from a regular graph with degree $s=q$ at most $r=5$ edges directed out of each vertex.
Thus we can invoke Theorem 2.2 to prove the claim of the corollary for $q$ odd.

If $q$ is even, we must also include self-loops to obtain a regular graph with degree $s=q+1$, from which we must now delete at most $r=6$ edges directed out of each vertex.
We again invoke Theorem 2.2 to prove the claim of the corollary for $q$ even.
\QED

Finally, we should point out that the regularity of these tessellations is unimportant.
McCann [M] has shown that cellular automata using majority voting and based on ``nice'' graphs
tolerate adversarial combined faults if the condition of Corollary 2.4 is satisfied by some even $q$
that is merely a lower bound to the degree of each vertex.
(A simple undirected graph is ``nice'' if it is connected, locally-finite, and discretely embeddable in the
plane.)
\sk

\heading{3. A Negative Result}

In this section, we shall obtain a lower bound $\Me(1/\xi^2)$
to the degree necessary to achieve fault tolerance
with fault rate $1/2 - \xi$.
The lower bound will be presented for binary automata, but the generalization
to more than two states is straightforward.
We do not assume monotonicity or self-duality of the transition functions.
Furthermore, our result applies even if the only faults are transient faults, and if they occur 
(independently in time and space) with probability {\it exactly\/} $\ep = 1/2 - \xi$ (that is,
when there is no adversary).

Our result depends on a lemma due to Evans and Schulman [E1, E2] that quantifies information loss
in circuits in which each gate fails independently with probability exactly $\ep = 1/2 - \xi$.
Their result is the culmination of a line work begun by Pippenger [P1] with a result applying to 
formulas (circuits in which each gate has ``fan-out'' one, so that the circuit forms a tree).
Pippenger's result was generalized to circuits by Feder [F2], and Feder's
result was quantitatively improved by Evans and Schulman.

If $X$ is a random variable taking values in a finite set $\calX$, we define the {\it entropy\/} $H(X)$ of $X$ by
$$H(X) = -\sum_{x\in\calX} \Pr[X=x]\log_2 \Pr[X=x].$$
We have $H(X)\le\log_2\#\calX$, with equality for and only for the uniform distribution.
If $X$ and $Y$ are random variables, we define their {\it mutual information\/} $I(X; Y)$ by
$$I(X; Y) = H(X) + H(Y) - H(X,Y).$$
We have $I(X; Y)\ge 0$ from the {\it subadditivity\/} $H(X,Y)\le H(X) + H(Y)$ of entropy.
For $0\le p\le 1$, we define $h(p) = -p\log_2 p - (1-p)\log_2 (1-p)$, the entropy of a random variable
that assumes the value $1$ with probability $p$ and the value $0$ with probability $1-p$.

\label{Lemma 3.1:}
(Evans and Schulman [E1, E2])
Consider a circuit with one input $a$ and one output $b$, 
and in which the output of every gate is complemented with
probability exactly $\ep$.
Let $a$ be fed by a random variable $X$ uniformly distributed on $\{0,1\}$, and let
$Y$ be the resulting random variable produced at $b$.
Then
$$I(X;Y) \le \sum_p (1-2\ep)^{2\abs{p}},$$
where the sum is over all paths $p$ from $a$ to $b$, and $\abs{p}$ denotes the length
(number of gates on) the path $p$.

Of crucial importance to us is the factor $2$ appearing in the exponent; it is exactly this factor
by which Evans and Schulman's result improves Feder's.

\label{Theorem 4.2:}
Let $M$ be a cellular automaton in which the transition function for each cell depends on the states
of at most $d$ neighbors.
Then if $M$ tolerates pure transient faults with fault rate $\ep = 1/2 - \xi$ and error probability at most
$\de < 1/2$, we must have 
$$d\ge 1/4\xi^2.$$

For the proof, we shall need the following special case of Fano's lemma.

\label{Lemma 4.3:}
(R.~M. Fano; see Fano [F1], \S 6.2)
Let $X$ and $Y$ be binary random variables, with $X$ uniformly distributed on $\{0,1\}$.
If $\Pr(X\not=Y)\le \de < 1/2$, then $I(X;Y)\ge 1 - h(\de)$.

\label{Proof:}
If $X\oplus Y$ denotes the exclusive-OR (sum modulo $2$) of $X$ and $Y$, then $X\oplus Y = 1$
if and only if $X\not=Y$.
We then have
$$\eqalign{
I(X; Y) &= H(X) + H(Y) - H(X,Y) \cr
&= 1 + H(Y) - H(X,Y) \cr
&= 1 + H(Y) - H(X\oplus Y,Y) \cr
&\ge 1 - H(X\oplus Y) \cr
&\ge 1 - h(\de). \cr
}$$
Here we have used the definition of $I(X;Y)$, the fact that $H(X) = 1$
(since $X$ is uniformly distributed on $\{0,1\}$), the identity
$H(X,Y) = H(X\oplus Y)$ (since any two of $X$, $Y$ and $X\oplus Y$ determine the third),
the subadditivity of entropy $H(X\oplus Y)\le H(X\oplus Y) + H(Y)$,
and the inequality $H(X\oplus Y)\le h(\de)$
(since $h(\de)$ is a non-decreasing function of $\de$ for $0\le \de\le 1/2$,
and $\Pr(X\oplus Y) = \Pr(X\not=Y)\le \de < 1/2$).
\QED

\label{Proof of Theorem 4.2:}
Given a binary cellular automaton, a cell $v$ and a time $t\ge 1$, we construct a circuit as follows.
The circuit will have a single input $a$, a single output $b$ and $t$ layers of gates.
The gates in a given layer will correspond to a finite subset of the cells in the automaton.
The $t$-th layer will contain a single gate, corresponding to the cell $v$,
and this gate will feed the output $b$.
For $s = t-1, \ldots, 2, 1$, the gates in the $s$-th layer will correspond to the cells that are neighbors of cells
corresponding to gates in the $(s+1)$-st layer, and the gates in the $(s+1)$-st layer will be fed by the
appropriate gates in the $s$-th layer.
All gates in the first layer are fed from the input $a$.

Suppose now that the input $a$ is fed a random variable $X$ uniformly distributed in $\{0, 1\}$,
suppose that  the  gates suffer faults (that is, that their outputs are complemented) independently 
with probability exactly $\ep$,
and let $Y$ be the random variable produced at the output $b$.
Suppose further that the cellular automaton is started with all initial states equal to $X$,
that the cellular automaton suffers pure transient faults (that is, states are complemented, independently in time and space)  with probability exactly $\ep$.
Then the distribution of the state of cell $v$ at time $t$ is the same as that of $Y$.

In this circuit, there are at most $d^t$ paths from $a$ to $b$, so
$$I(X; Y)\le d^t\,(2\xi)^{2t}$$
by Lemma 3.1.
Since
$$I(X;Y)\ge 1 - h(\de)$$
by Lemma 3.3, we obtain
$$d \ge \(1 - h(\de)\)^{1/t}\big/4\xi^2.$$
Since $\de<1/2$, we have $1 - h(\de) > 0$, so
$\(1 - h(\de)\)^{1/t}\to 1$ as $t\to\infty$.
Thus we obtain the desired bound
$$d\ge 1/4\xi^2.$$
\QED
\sk

\heading{4. Conclusion}

We have obtained nearly matching upper and lower bounds on the degree required 
by cellular automata to tolerate fault rates approaching $1/2$.
We have confined our attention to the binary case, but all of our results generalize easily
to the case of an arbitrary finite set of states.

Two questions are left unanswered by this work.
The first, of course , concerns the logarithmic gap between the upper bound
$O\((1/\xi^2)\log(1/\xi)\)$ and the lower bound $O(1/\xi^2)$.
The second arises from the fact that our upper bounds apply only to automata based
on graphs that contain, in an appropriate sense, infinite regular trees.
These graphs have natural embeddings in the hyperbolic plane.
We do not know whether fault rates approaching $1/2$ can be tolerated automata
in Euclidean spaces, even with dimensions higher than two or three.
The known fault-tolerance results for automata in Euclidean spaces
(see Toom [T3], for example) require that the fault rate be ``sufficiently small'',
and the fault-rate threshold is not decreased by increasing the degree.

Finally, we should point out that in our upper-bound results, for trees and other regular tessellations
of the hyperbolic plane, we have not considered
any transition functions other than those based on majority voting among all neighbors,
which is symmetric under all automorphisms of the underlying graph.
It is known, however, that in other contexts (see Pippenger [P3])
asymmetric transition functions are able
to achieve fault tolerance when symmetric functions cannot.
\sk

\heading{5. Acknowledgment}

The research reported here was supported
by Grant CCF 0430656 from the National Science Foundation.
\sk

\heading{6.  References}

\refbook B; E. R. Berlekamp, J. H. Conway and R. K. Guy;
Winning Ways for Your Mathematical Plays; Academic Press, 1982, v.~2.

\item{[C1]} H. S. M. Coxeter,
``Regular Honeycombs in Hyperbolic Space'', in:
{\it Proceedings of the International Congress of Mathematicians, 1954},
North-Holland Publishing, 1956, v.~III, pp.~155--169
(reprinted in H.~S.~M. Coxeter, {\it The Beauty of Geometry},  Dover Publications, 1999).

\refbook C2; H. S. M. Coxeter and W. O. J. Moser;
Generators and Relations for Discrete Groups;
4th edition, Springer-Verlag, 1980.

\refbook E1; W. S. Evans;
Information Theory and Noisy Computation;
Ph.\ D. Thesis, Department of Computer~Science,
University of California at Berkeley, 1994.

\ref E2;  W. S. Evans and L. J. Schulman;
``Signal Propagation and Noisy Circuits'';
IEEE Trans.\ Inform.\ Theory; 45:7 (1999) 1--7.

\refbook F1; R. M. Fano;
Transmission of Information;
MIT Press, 1961.

\ref F2; T. Feder;
``Reliable Computation by Networks in the Presence of Noise'';
IEEE Trans.\ Inform.\ Theory; 35:3 (1989) 569--571.

\refinbook G; P. G\'{a}cs;
``Self-Correcting Two-Dimensioal Arrays'';
in: S.~Micali (Ed.);
Randomness and Computation;
JAI Press, 1989, v.~5, pp.~223--326.

\ref K; A. V. Kuznetsov;
``Information Storage in a Memory Assembled from Unreliable Components'';
Problems of Information Transmission; 9:3 (1973) 254--264
 (translated from {\it Problemy Peredachi Informatsi\u{\i}}, 9:3 (1973) 100--114).

\refbook M; M. A. McCann;
Memory in Media with Manufacturing Faults;
Ph.~D. Thesis, Department of Computer Science,
Princeton University, September 2007.

\refinbook N1; J. von Neumann;
``Probabilistic Logics and the Synthesis of Reliable Organisms from Unreliable Components'';
in: C.~E. Shannon and J.~McCarthy (Ed's);
Automata Studies; Princeton University Press, 1956, pp.~43--98.

\refbook N2; J. von Neumann (compiled by A.~W. Burks);
 Theory of Self-Reproducing Automata;
 University of Illinois Press, 1966.

\ref P1; N. Pippenger;
``Reliable Computation by Formulas in the Presence of Noise'';
IEEE Trans.\ Inform.\ Theory; 34:2 (1988) 194--197.

\ref P2; N. Pippenger;
``Invariance of Complexity Measures for Networks with Unreliale Gates'';
J. Assoc.\ Comput.\ Mach.; 36:3 (1989) 531--539.

\ref P3; N. Pippenger;
``Symmetry in Self-Correcting Cellular Automata'';
Journal of Computer and System Sciences; 49:1 (194) 83--95.

\ref S; O. N. Stavskaya and I. I. Pyatetski\u{\i}-Shapiro;
``On Homogeneous Nets of Spontaneously Active Elements'';
Systems Theory Research; 20 (1976) 75--88
 (translation of {\it Problemy Kibernetiki}, 20 (1968) 91--106).
 
 \ref T1; A. L. Toom;
 ``Nonergodic Multidimensional Systems of Automata'';
 Problems of Information Transmission; 10:3 (1974) 239--246
 (translated from {\it Problemy Peredachi Informatsi\u{\i}}, 10:3 (1974) 70--79).
 
  \ref T2; A. L. Toom;
 ``Monotonic Binary Cellular Automata'';
 Problems of Information Transmission; 12:1 (1976) 33--37
 (translated from {\it Problemy Peredachi Informatsi\u{\i}}, 12:1 (1976) 48--54).
 
 \refinbook T3; A. L. Toom;
 ``Stable and Attractive Trajectories in Multicomponent Systems'';
 in: R.~L. Dobrushin and Ya.~A. Sinai (Ed's);
 Multicomponent Random Systems; Marcel Dekker, 1980, pp.~549--575.
 
 \refinbook U; S. Ulam;
 ``Random Processes and Transformations'';
 in: L.~M. Graves, E.~Hille, P.~Smith and O.~Zariski (Ed's);
 Proceedings of the International Congress of Mathematicians, 1950;
 American Mathematical Society, 1952, v.~2, pp.~264--275.

\bye